\documentclass[11pt]{amsart}
\usepackage{amssymb,amsfonts}
\usepackage{euscript,mathrsfs}
\usepackage{latexsym}
\usepackage{xspace}
\usepackage{amscd}
\usepackage{amsmath}
\usepackage{color}
\usepackage{cite}
\usepackage{graphicx}
\usepackage[utf8]{inputenc}
\usepackage{amsmath,amsthm,amsopn,graphicx,microtype}
\usepackage{url}
\calclayout

\theoremstyle{plain}

\theoremstyle{definition}

\theoremstyle{remark}

\newcommand{\Mac}{\emph{Macaulay2}\xspace}
\newcommand{\SOS}{\textsc{SumsOfSquares}\xspace}
\newcommand{\SDP}{\textsc{SemidefiniteProgramming}\xspace}

\newcommand{\QQ}{\mathbb{Q}}

\newcommand{\RR}{\mathbb{R}}
\newcommand{\kk}{\mathbb{K}}

\begin{document}

\title[Sums of squares in Macaulay2]{Sums of squares in Macaulay2}

\author{Diego Cifuentes}
\address{Massachusetts Institute of Technology \\ Cambridge, MA, USA}
\email{diegcif@mit.edu}

\author{Thomas Kahle}
\address{Otto-von-Guericke University \\ Magdeburg, Germany}
\email{thomas.kahle@ovgu.de}

\author{Pablo A. Parrilo}
\address{Massachusetts Institute of Technology \\ Cambridge, MA, USA}
\email{parrilo@mit.edu}

\begin{abstract}
  The \Mac package \SOS decomposes polynomials as sums of squares.
  It is based on methods to rationalize sum-of-squares decompositions due to Parrilo and Peyrl.
  The package features a data type for sums-of-squares polynomials, support for external semidefinite programming solvers, and optimization over varieties.
\end{abstract}

\maketitle

\section{Introduction}
\label{s:intro}

Let $\kk \!=\! \QQ$ or $\kk \!=\! \RR$ be the rational or real numbers and $R = \kk[x_{1},\dots,x_{n}]$ be the polynomial ring.
An element $f\!\in\! R$ is \emph{nonnegative} if $f(x) \!\ge\! 0$ for all $x \!\in\! \RR^{n}$,
and $f$ is a \emph{sum of squares} (SOS) if there are polynomials $f_{1},\dots,f_{m} \!\in\! R$ and positive scalars $\lambda_{1},\dots,\lambda_{m} \!\in\! \kk$ such that $f=\sum_{i}\lambda_i f_{i}^{2}$.
The scalars are not necessary when the field is $\kk\!=\!\RR$.
Clearly, a sum of squares is nonnegative, but not every nonnegative polynomial is a sum of squares.
Hilbert showed that the nonnegative polynomials of degree $d$ in $n$ variables are sums of squares if and only if: $n=1$; or $d=2$; or $n=2$ and $d=4$.
For an introduction to the area we recommend~\cite{scheiderer2009positivity,blekherman2012semidefinite}.

The \SOS package contains methods to compute sums-of-squares in \Mac~\cite{macaulay2}.
A particular focus is on trying to find rational sums-of-squares decompositions of polynomials with rational coefficients (whenever they exist).

Consider the basic problem of deciding whether a polynomial is a sum of squares.
Let $f \!\in\! R$ of degree~$2d$,
and $v \!\in\! R^N$ a vector whose entries are the $N \!=\! \binom{n+d}{d}$ monomials of degree $\leq\!d$.
The following fundamental result holds:
\begin{align*}
  f \text{ is SOS } 
  \qquad\iff\qquad
  \exists\, Q\in \mathbb{S}_+^N \text{ such that } f = v^T Q v,
\end{align*}
where $\mathbb{S}_+^N$ is the cone of $N{\times} N$ symmetric positive
semidefinite matrices; see~\cite[\S3.1]{blekherman2012semidefinite}.
This reduces the problem to finding a \emph{Gram matrix} $Q$ as above,
which can be done efficiently with \emph{semidefinite programming}
(SDP).

The method \verb|solveSOS| performs the computation above.
We use it here to verify that $f = 2 x^4 {+} 5 y^4 {-} 2 x^2 y^2 {+} 2 x^3 y$ is a sum of squares:
{\small
\begin{verbatim}
i1 : R = QQ[x,y];
i2 : f = 2*x^4+5*y^4-2*x^2*y^2+2*x^3*y;
i3 : sol = solveSOS f;
Executing CSDP
Status: SDP solved, primal-dual feasible
\end{verbatim}
}

\noindent
The ``Status'' line indicates that a Gram matrix was found, so $f$ is
indeed a sum of squares.  In the example above the package called an
external program to serve as semidefinite programming solver.  The
default solver is the open source program
CSDP~\cite{borchers1999csdp}, which is included in \Mac.  The output
of \verb|solveSOS| is an object of type \verb|SDPResult|.  It
contains, in particular, the Gram matrix~$Q$ and the monomial
vector~$v$.  {\small
\begin{verbatim}
i4 : (Q,v) = ( sol#GramMatrix, sol#Monomials )
o4 = ( | 2      1     -83/40 |,   | x2 | )
       | 1      43/20 0      |    | xy | 
       | -83/40 0     5      |    | y2 | 
\end{verbatim}
}
\noindent
The result of the semidefinite programming solver is a floating point approximation of the Gram matrix.
The \SOS package attempts to find a close enough rational Gram matrix by rounding its entries~\cite{peyrl2008computing}.
If this rounding procedure fails to find a feasible rational matrix, the method returns the floating point solution.
The procedure is guaranteed to work when the floating point Gram matrix lies in the interior of~$\mathbb{S}_+^N$.
See Appendix~\ref{sec:rounding} for more details about rational rounding.

The method \verb|sosPoly| extracts the sum-of-squares decomposition from the returned \verb|SDPResult|.
This is done via an LDL factorization (a variant of Cholesky factorization) of the Gram matrix.
For the function $f$ from above we get three squares:
{\small
\begin{verbatim}
i5 : s = sosPoly sol
            83 2    2 2    43  20 2       2    231773   2 2
o5 = (5)(- ---x  + y )  + (--)(--x  + x*y)  + (------)(x )
           200             20  43              344000
\end{verbatim}
}
\noindent
The output above is an object of type \verb|SOSPoly|.
An object of this type stores the coefficients $\lambda_{i}$ and polynomials (or generators) $f_{i}$  such that $f = \sum_{i}\lambda_{i}f_{i}^{2}$.
We can extract the coefficients and generators as follows:
{\small
\begin{verbatim}
i5 : coefficients s
o5 = {5, 43/20, 231773/344000}
i6 : gens s
o6 = {-83/200*x^2 + y^2, 20/43*x^2 + x*y, x^2}
\end{verbatim}
}

The method \verb|solveSOS| can also compute sums-of-squares decompositions in quotient rings.
This can be useful to prove nonnegativity of a polynomial on a variety.
We take an example from~\cite{parrilo2005exploiting}.
Consider proving that $f = 10{-}x^2{-}y$ is nonnegative on the circle defined by $g = x^2 {+} y^2 {-} 1$.
To do this, we check if $f$ is a sum of squares in the quotient ring $\QQ[x,y]/\langle g\rangle$.
For such a computation, an even degree bound must be given by the user, 
as otherwise it is not obvious how to choose the monomial vector~$v$.
In the following example we use $2d=2$ as the degree bound.
{\small
\begin{verbatim}
i1 : R = QQ[x,y]/ideal(x^2 + y^2 - 1);
i2 : f = 10-x^2-y;  d = 1;
i3 : sosPoly solveSOS (f, 2*d, TraceObj=>true)
Executing CSDP
Status: SDP solved, primal-dual feasible
            1       2    35    2
o3 = (9)(- -- y + 1)  + (--)(y)
           18            36
\end{verbatim}
}

\noindent
In the computation above the option \verb|TraceObj=>true| was used to reduce the number of squares in the SOS decomposition (see Section~\ref{s:arguments}).

\section{Sums of squares in ideals}
Let $I \subset \kk[x_{1},\dots,x_{n}]$ be an ideal.  
Given an even bound $2d$, consider the problem of finding a nonzero sum-of-squares polynomial of degree $\leq \!2d$ in the ideal~$I$.
If one of the generators of $I$ has degree $\leq \!d$, then the problem is trivial.
But otherwise the problem can be hard.
The method \verb|sosInIdeal| can be used to solve it.
One of the main motivations for this problem is that it reveals information about the \emph{real radical} of the ideal~$I$, 
i.e., the vanishing ideal of the real zeros of~$I$.
Indeed, if $f = \sum \lambda_i f_i^2 \in I$ then each of the factors $f_i$ must lie in the real radical of~$I$.

Given generators of the ideal $I=\langle h_1,\dots,h_m\rangle$, we may solve this problem by looking for some polynomial multipliers $l_i(x)$ such that  $\sum_i l_i(x) h_i(x)$ is a sum of squares.
The method \verb|sosInIdeal| can find these multipliers.
The input is a matrix containing the generators, and the degree bound~$2d$.
We illustrate this for the ideal 
$I=\langle x^2{-}4 x{+}2 y^2, 2 z^2{-}y^2{+}2 \rangle$
{\small
\begin{verbatim}
i1 : R = QQ[x,y,z];  d = 1;
i2 : h = matrix {{x^2-4*x+2*y^2, 2*z^2-y^2+2}};
i3 : (sol,mult) = sosInIdeal (h, 2*d);
i4 : sosPoly sol
      395    1       2    395    2
o4 = (---)(- - x + 1)  + (---)(z)
       2     2             2
i5 : h * mult == sosPoly sol
o5 = true
\end{verbatim}
}
\noindent
Another way to approach this problem is to construct the quotient
$S = \kk[x_{1},\dots,x_{n}]/I$ and then write $0\in S$ as a sum of
squares.  In this case the input to \verb|sosInIdeal| is simply the
quotient ring~$S$.  {\small
\begin{verbatim}
i6 : S = R/ideal h;
i7 : sosPoly sosInIdeal (S, 2*d);
      1031833    1       2    1031833    2
o7 = (-------)(- - x + 1)  + (-------)(z)
        2048     2              2048
\end{verbatim}
}
\noindent
In both cases we obtained a multiple of the sum-of-squares polynomial $(\frac{1}{2}x{-}1)^2{+}z^2$.
This computation reveals that $x{-}2,z$ lie in the real radical of~$I$.
Indeed, we have $\sqrt[\RR]{I} = \langle x{-}2,z,y^2{-}2\rangle$.

\section{SOS decompositions of ternary forms}

Hilbert showed that any nonnegative form $f\in \kk[x,y,z]$ can be decomposed as a quotient of sums of squares.
We can obtain this decomposition by iteratively calling \verb|sosInIdeal|.
Specifically, one can first find a multiplier $q_{1}$ such that $q_{1}f$ is a sum of squares.
Since $q_1$ is also nonnegative, we can then search for a multiplier $p_{1}$ such that $p_{1}q_{1}$ is a sum of squares, and so on.
The main observation is that the necessary degree of $p_{1}$ is lower than that of $q_{1}$~\cite{de2004products}.
Hence this procedure terminates, and we can write
\[
  f = \frac{p_{1}\cdots p_{s}}{q_{1}\cdots q_{t}} \qquad \text {
    $p_{i},q_{i}$ SOS}.
\]

As an illustration, we write the Motzkin polynomial as a quotient of sums of squares.
We first use the function \verb|library|, which contains a small library of interesting nonnegative forms.
{\small
\begin{verbatim}
i1 : R = QQ[x,y,z]
i2 : f = library ("Motzkin", {x,y,z})
      4 2    2 4     2 2 2    6
o2 = x y  + x y  - 3x y z  + z
\end{verbatim}
}
\noindent
We now apply the function \verb|sosdecTernary|, which implements the iterative algorithm from above.
{\small
\begin{verbatim}
i3 : (Nums,Dens) = sosdecTernary f;
Executing CSDP
i4 : num = first Nums
      2267   2 2    4 2    2003    1013 3     990   3        2 2    
o4 = (----)(x y  - z )  + (----)(- ----x y - ----x*y  + x*y*z )  + ... 
       64                   64     2003      2003                 
i5 : den = first Dens
      2267    2    1079    2    33    2
o5 = (----)(z)  + (----)(x)  + (--)(y)
       64           64           2
\end{verbatim}
}
\noindent
The result consists of two sums of squares, the second being the denominator.
We can check the computation as follows.
{\small
\begin{verbatim}
i6 : f*value(den) == value(num)
o6 = true
\end{verbatim}
}

\section{Parametric SOS problems}

The \SOS package can also solve parametric problems.
Assume now that $x \mapsto f(x;t)$ is a polynomial function, that depends affinely on some parameters~$t$.
The command \verb|solveSOS| can be used to search for values of the parameters such that the polynomial is a sum of squares.
In the following example, we change two coefficients of the Robinson polynomial so that it becomes a sum of squares.
{\small
\begin{verbatim}
i1 : R = QQ[x,y,z][s,t];
i2 : g = library("Robinson", {x,y,z}) + s*x^6 + t*y^6;
i3 : sol = solveSOS g;
Executing CSDP
Status: SDP solved, primal-dual feasible
i4 : sol#Parameters
o4 = | 34 |
     | 34 |
\end{verbatim}
}

\noindent
In the code above, the ring construction (first line) indicates that $s,t$ should be treated as parameters.
The values obtained were $s=t=34$.

It is also possible find the values of the parameters that optimize a given linear function.
This allows us to find lower bounds for a polynomial function $f(x)$,
by finding the largest $t$ such that $f(x)-t$ is a sum of squares.
Here we apply this method to the dehomogenized Motzkin polynomial.

{\small
\begin{verbatim}
i1 : R = QQ[x,z][t];
i2 : f = library ("Motzkin", {x,1,z});
i3 : sol = solveSOS (f-t, -t, RoundTol=>12);
Executing CSDP
Status: SDP solved, primal-dual feasible
i4 : sol#Parameters
o4 = | -729/4096 |
\end{verbatim}
}

\noindent
Alternatively, the method \verb|lowerBound| can be called with input~$f(x)$.
The method internally declares a new parameter~$t$ and optimizes $f(x)-t$.
{\small
\begin{verbatim}
i1 : R = QQ[x,z];
i2 : f = library ("Motzkin", {x,1,z});
i3 : (t,sol) = lowerBound (f, RoundTol=>12);
Executing CSDP
Status: SDP solved, primal-dual feasible
i4 : t
o4 = - 729/4096
\end{verbatim}
}

\section{Polynomial optimization}

In applications one often needs to find lower bounds for polynomials subject to some polynomial constraints.
More precisely, consider the problem
\begin{align*}
  \min_{x\in \RR^n} \quad f(x)
  \quad \text{ such that }\quad
  h_1(x)=\dots=h_m(x)=0,
\end{align*}
where $f, h_1,\dots,h_m$ are polynomials.
The \SOS package provides two ways to compute a lower bound for such a problem.
The most elegant approach is to construct the associated quotient ring, and then call \verb|lowerBound|.
This will look for the largest $t$ such that $f(x)-t$ is a sum of squares (in the quotient ring).
A degree bound~$2d$ must be given by the user.

{\small
\begin{verbatim}
i1 : R = QQ[x,y]/ideal(x^2 - x, y^2 - y);
i2 : f = x - y;  d = 1;
i3 : (t,sol) = lowerBound(f,2*d);
Executing CSDP
Status: SDP solved, primal-dual feasible
i4 : t
o4 = -1
i5 : f - t == sosPoly sol
o5 = true
\end{verbatim}
}

Calling \verb|lowerBound| as above is conceptually simple, but requires knowledge of a Gröbner basis, which is computed when constructing the quotient ring.
If no Gröbner basis is available there is an alternative way to call \verb|lowerBound| with just the equations $h_1,\dots,h_m$ as the input.
The method will then look for polynomial multipliers $l_i(x)$ such that $f(x) - t + \sum_i l_i(x)h_i(x)$ is a sum of squares.
This may result in larger semidefinite programs and weaker bounds.

{\small
\begin{verbatim}
i1 : R = QQ[x,y];
i2 : f = x - y;  d = 1;
i3 : h = matrix{{x^2 - x, y^2 - y}};
i4 : (t,sol,mult) = lowerBound (f, h, 2*d);
Executing CSDP
Status: SDP solved, primal-dual feasible
i5 : t
o5 = -1
i6 : f - t + h*mult == sosPoly sol
o6 = true
\end{verbatim}
}

Lower bounds for polynomial optimization problems critically depend on the degree bound chosen.
While higher degree bounds lead to better bounds, the computational complexity escalates quite rapidly.
Nonetheless, low degree SOS lower bounds often perform very well in applications.
In some cases, the minimizer might also be recovered from the \verb|SDPResult| with the method \verb|recoverSolution|.

{\small
\begin{verbatim}
i7 : recoverSolution sol
o7 = {x => 1.77345e-9, y => 1}
\end{verbatim}
}

\section{Optional arguments}
\label{s:arguments}

\subsection*{SDP Solver}
The optional argument \verb|Solver| is available for many package methods and a particular semidefinite programming solver can be picked by setting it.
These solvers are interfaced via the auxiliary \Mac package \SDP~\cite{sdpM2}.
The package provides interfaces to the open source solvers CSDP~\cite{borchers1999csdp} and SDPA~\cite{yamashita2003implementation}, and the commercial solver MOSEK~\cite{mosek}.
There is also a built-in solver in the \Mac language.
In our experience CSDP and MOSEK give the best results.
CSDP is provided as part of \Mac and configured as the default.

\subsection*{Rounding tolerance}
The method \verb|lowerBound| has the optional argument
\verb|RoundTol|, which specifies the precision of the rational
rounding.  Smaller values of \verb|RoundTol| lead to rational matrices
with smaller denominators but farther from the numerical solution.
The rational rounding may be skipped by setting it to infinity.

\subsection*{Trace objective}
The option \verb|TraceObj| tells the semidefinite programming solver to minimize the trace of the Gram matrix.
This is a known heuristic to reduce the number of squares in the SOS decomposition.

\appendix

\section{Rational rounding}%
\label{sec:rounding}

Sums-of-squares problems are solved numerically using an semidefinite programming solver,
and afterwards the package attempts to round the floating point solution to rational numbers.
We briefly describe the rounding procedure, which was proposed in~\cite{peyrl2008computing}.

Let $f\in \QQ[x_1,\dots,x_n]$ and consider the affine space
$\mathcal{L}:= \{Q: v^T\! Q v \!=\! f\}$, where $v$ is a given
monomial vector.  A Gram matrix is an element of
$\mathcal{L}\cap \mathbb{S}_+^N$.  The semidefinite programming solver
returns a numerical matrix~$Q_n$, an ``approximate'' Gram matrix,
which may not lie exactly on~$\mathcal{L}$.  The rounding problem
consists in finding a nearby Gram matrix~$Q_r$ with rational entries.

The procedure from~\cite{peyrl2008computing} consists of two steps.
First, the entries of~$Q_n$ are rounded to a rational matrix~$Q_r'$.
Then $Q_r$ is obtained as the orthogonal projection of $Q_r'$ onto~$\mathcal{L}$.
The image of the projection is rational, lies in~$\mathcal{L}$, but need not be positive semidefinite.
We may ensure that $Q_r\in \mathbb{S}_+^N$ if the numerical matrix $Q_n$ is in the interior of~$\mathbb{S}_+^N$ and sufficiently close to $\mathcal{L}$.
More precisely, assume that $\lambda$, the smallest eigenvalue of~$Q_n$, is greater than the distance $\delta:= \textrm{dist}(Q_n,\mathcal{L})$.
Then setting the {rounding tolerance} $\textrm{dist}(Q_n,Q_r')$ smaller than $\sqrt{\lambda^2 - \delta^2}$ guarantees that $Q_r \in \mathbb{S}_+^N$; 
see \cite[Prop.~8]{peyrl2008computing}.
\section*{Acknowledgment}
\label{sec:acknowledgement}
The authors would like to thank Bernd Sturmfels and the Max-Planck-Institute für
Mathematik in den Naturwissenschaften in Leipzig for hosting the \Mac workshop in May 2018.
We thank Ilir Dema, Nidhi Kainsa and Anton Leykin, who contributed to the code.
The package code contains parts of a proof-of-concept implementation of the methods in \cite{peyrl2008computing}, originally written by Pablo Parrilo and Helfried Peyrl.
Parts of this work were done while Diego Cifuentes visited the Max-Planck-Institute MiS, and while the first two authors visited ICERM supported by NSF grant No. DMS-1439786.
Thomas Kahle is supported by the German Research Foundation under grant 314838170, GRK 2297 MathCoRe.
Pablo Parrilo is supported in part by the National Science Foundation through NSF Grant CCF-1565235.

\bibliographystyle{amsplain}
\bibliography{sos}

\end{document}